\documentclass[12pt]{article}
\usepackage[cp866]{inputenc}
\usepackage{citehack}
\usepackage{amsmath}
\usepackage{amsfonts}
\usepackage{amssymb}
\usepackage{latexsym}
\newtheorem{theorem}{Theorem}

\newtheorem{defin}{Definition}
\newtheorem{lem}{Lemma}

\def\i{{\hat{i}}}
\def\j{{\hat{j}}}
\def\k{{\hat{k}}}
\def\l{{\hat{l}}}
\def\ii{{\hat{\hat{i}}}}
\def\jj{{\hat{\hat{j}}}}
\def\kk{{\hat{\hat{k}}}}
\def\ll{{\hat{\hat{l}}}}
\def\iv{{\tilde{i}}}
\def\jv{{\tilde{j}}}
\def\kv{{\tilde{k}}}
\def\lv{{\tilde{l}}}
\def\iiv{{\tilde{\tilde{i}}}}
\def\jjv{{\tilde{\tilde{j}}}}
\def\kkv{{\tilde{\tilde{k}}}}
\def\llv{{\tilde{\tilde{l}}}}

\def\zr{\ltimes}

\def\Real{\mathbb{R}}
\def\g{\mathfrak{g}}
\def\h{\mathfrak{h}}
\def\hol{\mathfrak{hol}}
\def\so{\mathfrak{so}}
\def\spin{\mathfrak{spin}}

\def\gl{\mathfrak{gl}}
\def\su{\mathfrak{su}}
\def\u{\mathfrak{u}}

\def\f{\mathfrak{f}}

\def\z{\mathfrak{z}}
\def\R{\mathcal{R}}
\def\P{\mathcal{P}}

\def\A{\mathcal {A}}

\def\N{\mathcal {N}}
\def\K{\mathcal {K}}
\def\id{\mathop\text{\rm id}\nolimits}
\def\spa{\mathop\text{{\rm span}}\nolimits}
\def\Hom{\mathop\text{\rm Hom}\nolimits}

\def\pr{\mathop\text{\rm pr}\nolimits}
\def\Ga{\Gamma}

\title{Metrics that realize all Lorentzian holonomy algebras}

\author{Anton S. Galaev \thanks{EMail: \tt  galaev@mathematik.hu-berlin.de}}
\begin{document}

\maketitle \rightline{Dedicated to Dmitri Vladimirovich Alekseevsky on his 65th birthday} \vskip-50ex

{\renewcommand{\abstractname}{Abstract}\begin{abstract}
All candidates to the weakly-irreducible not irreducible holonomy algebras of Lorentzian manifolds are known.
In the present paper metrics that realize all these candidates as holonomy algebras are given.
This completes the classification of the Lorentzian holonomy algebras. Also new examples of metrics with the holonomy algebras
$g_2\zr\Real^7\subset\so(1,8)$ and $\spin(7)\zr\Real^8\subset\so(1,9)$ are constructed.
\end{abstract}


{\bf Mathematical subject codes:} 53C29, 53C50, 53B30

\section*{Introduction}

The classification of the holonomy algebras for Riemannian manifolds is a well-known classical result.
By the  Borel-Lichnerowicz theorem a Riemannian manifold is locally a product of
Riemannian manifolds with irreducible holonomy algebras (\cite{Bo-Li}).
In 1955 M. Berger gave a list of  possible irreducible holonomy algebras of Riemannian manifolds, see \cite{Ber}.
Later, in 1987 R. Bryant constructed  metrics for the  exceptional algebras of this list, see \cite{Bryant}.
For more details see \cite{Al,Besse,Joyce}.

The classification problem for the  holonomy algebras of pseudo-Riemannian manifolds is still open.
The main difficulty is that the holonomy algebra can preserve an isotropic subspace of the tangent space.
A subalgebra $\g\subset\so (r,s)$ is called {\it weakly-irreducible} if it does not  preserve any nondegenerate
proper subspace of $\Real^{r,s}$.
The Wu theorem states that {\it a pseudo-Riemannian manifold is locally a product of
pseudo-Riemannian manifolds with weakly-irreducible holonomy algebras}, see \cite{Wu}.
So, it is enough to consider only weakly-irreducible holonomy algebras.
If a holonomy algebra is irreducible, then it is weakly-irreducible. In \cite{Ber}
M. Berger gave also  a classification of possible  irreducible holonomy algebras for pseudo-Riemannian manifolds,
but there is no classification for weakly-irreducible not irreducible holonomy algebras for general pseudo-Riemannian manifolds.
About the Lorentzian case see below. There are some partial results for holonomy algebras of pseudo-Riemannian manifolds of
signature $(2,N)$, see \cite{Ik2n,G3,G4}, and $(N,N)$, see \cite{B-Inn}.

We consider the  holonomy algebras of Lorentzian manifolds. From Berger's list it follows that {\it the
only irreducible holonomy algebra of Lorentzian
manifolds is  $\so(1,n+1)$}, see \cite{Disc-Ol} and \cite{Bo-Ze} for direct proofs of this fact.
In 1993  L. Berard Bergery and A. Ikemakhen classified  weakly-irreducible not irreducible
subalgebras of $\so(1,n+1)$, see \cite{B-I}. More precisely, they divided these subalgebras into 4 types, and associated to each
such subalgebra $\g\subset\so(1,n+1)$ a subalgebra $\h\subset\so(n)$, which is called the orthogonal part of $\g$.
In \cite{G2} more geometrical proof of this result was given.
The Lie algebras of type 1 and 2 have the forms $(\Real\oplus\h)\zr\Real^n$ and $\h\zr\Real^n$ respectively. The Lie algebras
of type 3 and 4 can be obtained from the first two by some twistings.
Recently T. Leistner proved that {\it if $\g\subset\so(1,n+1)$ is the holonomy
algebra of a Lorentzian manifold, then its orthogonal part $\h\subset\so(n)$ is the holonomy algebra of a Riemannian manifold},
see \cite{Le1,Le2,Le3,Le4}.
This gives a list of possible holonomy algebras for Lorentzian manifolds (Berger algebras). To complete the classification
of holonomy algebras one must prove that all Berger algebras can be realized as holonomy algebras.  In \cite{B-I} were
given metrics that realize all Berger algebras of type 1 and 2.
In \cite{BoubelPhD,Boubel} Ch. Boubel studied possible shapes of local metrics
for Lorentzian manifolds with weakly-irreducible not irreducible holonomy algebras.
In particular, he gave equivalent conditions for such manifolds
to have the holonomy of type 1,2,3 or 4 and  parameterized the set of
germs of metrics giving a holonomy algebra of each type.
In \cite{LS5} K. Sfetsos and D. Zoakos constructed metrics with the holonomy algebras
$\su(2)\zr\Real^4\subset\so(1,5)$, $\su(3)\zr\Real^6\subset\so(1,7)$ and $g_2\zr\Real^7\subset\so(1,8)$.

In the present paper we construct metrics that realize all Berger algebras as holonomy algebras. The method of the construction generalizes
an example of A. Ikemakhen given in \cite{Ik}. The coefficients of the constructed metrics are polynomial functions, hence the holonomy algebra at a point
is generated by the values of the curvature tensor and of its derivatives at this point, and it can be computed.
This completes the classification of the  holonomy algebras of Lorentzian manifolds.
As application we construct new examples of  metrics with the holonomy algebras
$g_2\zr\Real^7\subset\so(1,8)$ and $\spin(7)\zr\Real^8\subset\so(1,9)$.

{\bf Acknowledgements.}
I would like to thank Dmitri Vladimirovich Alekseevsky for introducing me to the  Lorentzian holonomy algebras
and for his careful attention to my work during the last four years. I am grateful to Charles Boubel, who
took my attention to the problem of construction of metrics. I thank the Erwin Schr\"odinger Institute, where
the work on this paper was finished.

\section{Preliminaries}

Let $(\Real^{1,n+1},\eta)$ be a Minkowski space of dimension $n+2$,
where  $\eta$ is a metric on $\Real^{n+2}$ of
signature $(1,n+1)$. We fix a basis
$p,e_1,...,e_n,q$ of $\Real^{1,n+1}$ such that the Gram matrix
of $\eta$ has the form $\left(\begin{array}{ccc} 0 & 0 & 1\\ 0 & E_n & 0 \\ 1 & 0 & 0 \\ \end{array}\right)$,
where $E_n$ is the $n$-dimensional identity matrix. We will denote by $\Real^n\subset\Real^{1,n+1}$ the Euclidean subspace spanned by the vectors $e_1,...,e_n$.

\begin{defin} A subalgebra $\g\subset \so(1,n+1)$ is
called  irreducible if it does not preserve any proper
subspace of $\Real^{1,n+1}$; $\g$ is called  weakly-irreducible if
it does not  preserve any nondegenerate proper subspace of $\Real^{1,n+1}$.
\end{defin}

Obviously, if $\g\subset \so(1,n+1)$ is irreducible, then it is weakly-irreducible.

Denote  by $\so(1,n+1)_{\Real p}$ the subalgebra of $\so(1,n+1)$ that preserves the isotropic line $\Real p$.
The Lie algebra $\so(1,n+1)_{\Real p}$ can be identified with the following matrix algebra
$$\so(1,n+1)_{\Real p}=\left\{\left. \left (\begin{array}{ccc} a &X & 0\\ 0 & A &-X^t \\ 0 & 0 & -a \\
\end{array}\right)\right|\, a\in \Real,\, X\in \Real^n,\,A \in \so(n) \right\} .$$
The above matrix can be identified with the triple $(a,A,X)$.
Define the following subalgebras of $\so(1,n+1)_{\Real p}$,  $\A=\{(a,0,0)|a\in \Real\},$ $\K=\{(0,A,0)|A\in \mathfrak{so}(n)\}$ and
$\mathcal N=\{(0,0,X)|X\in \Real^n\}$.
We see that $\A$ commutes with $\K$, and $\N$ is a commutative ideal. We also see that $$[(a,A,0),(0,0,X)]=(0,0,aX+AX).$$
We have the decomposition $$\so(1,n+1)_{\Real p}=(\mathcal A\oplus\mathcal K)\zr\mathcal N=(\Real\oplus\so(n))\zr\Real^n.$$

If a weakly-irreducible subalgebra $\g\subset \so(1,n+1)$ preserves a degenerate proper subspace $U\subset \Real^{1,n+1}$, then it preserves the
isotropic line $U\cap U^\bot$, and $\g$ is conjugated to a weakly-irreducible subalgebra of $\so(1,n+1)_{\Real p}$.

Let $\h\subset\so(n)$ be a subalgebra. Recall that $\h$ is a compact Lie
algebra and we have the decomposition $\h=\h'\oplus\z(\h)$, where $\h'$ is the commutant of $\h$ and $\z(\h)$ is the center of $\h$ (\cite{V-O}).

The following result is due to L. Berard Bergery and A. Ikemakhen.

{\bf Theorem \cite{B-I}}\label{B-I}
{\it A subalgebra $\g\subset \so(1,n+1)_{\Real p}$ is  weakly-irreducible if and only if
$\g$ belongs to one of the following types
\begin{description}
\item[type 1.] $\g^{1,\h}=(\Real\oplus\h)\zr\Real^n=\left\{\left. \left (\begin{array}{ccc}
a &X & 0\\ 0 & A &-X^t \\ 0 & 0 & -a \\
\end{array}\right)\right|\, a\in \Real,\,\,A \in \h,\, X\in \Real^n
 \right\}$, where $\h\subset\so(n)$ is a subalgebra;

\item[type 2.] $\g^{2,\h}=\h\zr\Real^n=\left\{\left. \left (\begin{array}{ccc}
0 &X & 0\\ 0 & A &-X^t \\ 0 & 0 & 0 \\
\end{array}\right)\right|\,\,A \in \h,\,  X\in \Real^n  \right\}$;

\item[type 3.] $\g^{3,\h,\varphi}=\{(\varphi(A),A,0)|A\in\h\}\zr\Real^n=\left\{\left. \left (\begin{array}{ccc}
\varphi(A) &X & 0\\ 0 & A &-X^t \\ 0 & 0 & -\varphi(A) \\
\end{array}\right)\right|\,\,A \in \h,\,  X\in \Real^n \right\}$,
where $\h\subset\so(n)$ is a subalgebra with $\z(\h)\neq\{0\}$, and  $\varphi :\h\to\Real$ is a non-zero linear map with $\varphi|_{\h'}=0$;

\item[type 4.] $\g^{4,\h,m,\psi}=\{(0,A,X+\psi(A))|A\in\h,X\in\Real^m\}\\=\left\{\left. \left (\begin{array}{cccc}
0 &X&\psi(A) & 0\\ 0 & A&0 &-X^t \\ 0 & 0 & 0 &-\psi(A)^t \\
0&0&0&0\\ \end{array}\right)\right|\,\,A \in \h,\,  X\in \Real^{m} \right\}$,
where  $0<m<n$ is an integer, $\h\subset\so(m)$ is a subalgebra with $\dim\z(\h)\geq n-m$,
and $\psi:\h\to \Real^{n-m}$ is a surjective linear map with $\psi|_{\h'}=0$. \end{description}}

\begin{defin}
The subalgebra $\h\subset\so(n)$ associated to a
weakly-irreducible subalgebra $\g\subset \so(1,n+1)_{\Real p}$ in the above theorem  is called  the
orthogonal part of $\g$.
\end{defin}

Let $(M,g)$ be a Lorentzian manifold of dimension $n+2$  and $\g$ the holonomy algebra (that is the Lie algebra of the holonomy group)
at a point $x\in M$. By Wu's theorem (see \cite{Wu}) $(M,g)$ is locally indecomposable, i.e. is not locally a product
of two pseudo-Riemannian manifolds if and only if the holonomy algebra  $\g$ is weakly-irreducible. If the holonomy algebra $\g$
is irreducible,  then $\g =\so(T_xM,g_x)$ (\cite{Ber}).  So we may assume that it is  weakly-irreducible and not irreducible. Then $\g$
preserves an isotropic line $\ell \subset T_xM$. We can identify the tangent space $T_xM$ with $\Real^{1,n+1}$ such that $g_x$ corresponds to $\eta$ and
$\ell$ corresponds to the line $\Real p$. Then $\g$ is identified with a weakly-irreducible  subalgebra of $\so(1,n+1)_{\Real p}$.

Let $W$ be a vector space and $\f\subset \gl(W)$ a subalgebra.

\begin{defin}
The vector  space
$$\R(\f)=\{R\in\Hom(W\wedge W,\f)|R(u\wedge v)w+R(v\wedge w)u+R(w\wedge u)v=0 \text{ for all } u,v,w\in W\}$$
is called the space of curvature tensors of type $\f$.
Denote by $L(\R(\f))$ the vector subspace of $\f$ spanned by $R(u\wedge v)$ for all $R\in\R(\f),$ $u,v\in W,$
$$L(\R(\f))=\spa\{R(u\wedge v)|R\in\R(\f),\,u,v\in W\}.$$
A subalgebra $\f\subset\gl(W)$ is called a Berger algebra if $L(\R(\f))=\f$.
\end{defin}

From the Ambrose-Singer theorem (\cite{Am-Si}) it follows that {\it if $\g\subset\so(1,n+1)_{\Real p}$ is the holonomy algebra
of a Lorentzian manifold, then $\g$ is a Berger algebra}.

Let $\h\subset\so(n)$ be a subalgebra.

\begin{defin}
The vector  space
$$\P(\h)=\{P\in \Hom (\Real^n,\h)|\,\eta(P(u)v,w)+\eta(P(v)w,u)+\eta(P(w)u,v)=0\text{ for all }u,v,w\in \Real^n\}$$
is called the space of  weak-curvature tensors of type $\h$.
A subalgebra $\h\subset\so(n)$ is called a weak-Berger algebra if $L(\P(\h))=\h$, where
$$L(\P(\h))=\spa\{P(u)|P\in\P(\h),\,u\in \Real^n\}$$ is the vector subspace of $\h$ spanned by
$P(u)$ for all $P\in\P(\h)$ and $u\in \Real^n$.
\end{defin}

The following theorem was proved in \cite{G1}.

{\bf Theorem \cite{G1}} {\it A weakly-irreducible subalgebra $\g\subset\so(1,n+1)_{\Real p}$ is a Berger algebra if and only if
its orthogonal part $\h\subset\so(n)$ is a weak-Berger algebra.}

Recently T. Leistner proved the following theorem.

{\bf Theorem \cite{Le1,Le2,Le3,Le4}}
{\it A subalgebra $\h\subset\so(n)$ is a weak-Berger algebra if and only if $\h$ is a Berger algebra.}

Recall that from the classification of Riemannian holonomy algebras it follows that {\it a subalgebra $\h\subset\so(n)$ is a Berger algebra if and only if
$\h$ is the holonomy algebra of a Riemannian manifold.}

Thus {\it a subalgebra $\g\subset\so(1,n+1)$ is a  weakly-irreducible not irreducible Berger algebra if and only if $\g$ is conjugated to one
of the subalgebras $\g^{1,\h},\g^{2,\h},\g^{3,\h,\varphi},\g^{4,\h,m,\psi}\subset\so(1,n+1)_{\Real p}$, where $\h\subset\so(n)$ is the
holonomy algebra of a Riemannian manifold.}

To complete the classification of Lorentzian holonomy algebras
we must realize all weakly-irreducible Berger subalgebra of $\so(1,n+1)_{\Real p}$
as the holonomy algebras. There are some examples.

{\bf Example \cite{B-I}.} In 1993 L. Berard Bergery and A. Ikemakhen realized
the weakly-irreducible Berger subalgebra of $\so(1,n+1)_{\Real p}$ of type 1 and 2
as the holonomy algebras. They constructed the following metrics.
Let $\h\subset\so(n)$ be the holonomy algebra of a Riemannian manifold.
Let $x^{0},x^{1},...,x^{n},x^{n+1}$ be the standard coordinates on $\Real^{n+2}$,
$h$ be a metric on $\Real^n$ with the holonomy algebra $\h$, and $f(x^{0},...,x^{n+1})$ be  a function with
$\frac{\partial f}{\partial x^1}\neq 0,$...$,\frac{\partial f}{\partial x^n}\neq 0$.
If $\frac{\partial f}{\partial x^0}\neq 0$, then the holonomy algebra of the metric
$$g=2dx^0dx^{n+1}+h+f\cdot(dx^{n+1})^2$$ is $\g^{1,\h}$. If  $\frac{\partial f}{\partial x^0}=0$, then the holonomy algebra of the metric $g$ is $\g^{2,\h}$.

\vskip0.3cm

In the next section we will construct metrics that realize all weakly-irreducible Berger algebras.
We will use the space $\P(\h)$ and the fact that $\h=L(\P(\h))$.
The idea of the constructions is given by the following example of A. Ikemakhen.

{\bf Example \cite{Ik}.}
Let $x^{0},x^{1},...,x^{5},x^{6}$ be the standard coordinates on $\Real^{7}$. Consider the following metric
$$g=2dx^0dx^{6}+\sum^{5}_{i=1}(dx^i)^2+2\sum^{5}_{i=1}u^i dx^{i} dx^{6},$$
where $$\begin{array}{ll}u^1=-(x^3)^2-4(x^4)^2-(x^5)^2,& u^2=u^4=0,\\ u^3=-2\sqrt{3}x^2x^3-2x^4x^5,& u^5=2\sqrt{3}x^2x^5+2x^3x^4. \end{array}$$

The holonomy algebra  of this metric at the point $0$ is $\g^{2,\rho(\so(3))}\subset\so(1,6)$, where
$\rho:\so(3)\to\so(5)$ is the representation given by the highest irreducible
component of the representation $\otimes^2\id:\so(3)\to\otimes^2\so(3)$.
The image $\rho(\so(3))\subset\so(5)$ is spanned by the matrices
$$A_1=\left(\begin{smallmatrix}0&0&-1&0&0\\0&0&\sqrt{3}&0&0\\1&-\sqrt{3}&0&0&0\\0&0&0&0&-1\\0&0&0&1&0\end{smallmatrix}\right),
A_2=\left(\begin{smallmatrix}0&0&0&-4&0\\0&0&0&0&0\\0&0&0&0&-2\\4&0&0&0&0\\0&0&2&0&0\end{smallmatrix}\right),
A_3=\left(\begin{smallmatrix}0&0&0&0&-1\\0&0&0&0&-\sqrt{3}\\0&0&0&-1&-1\\0&0&1&0&0\\1&\sqrt{3}&1&0&0\end{smallmatrix}\right).$$

We have $\pr_{\so(n)}\left(R\left(\frac{\partial}{\partial x^{3}},\frac{\partial}{\partial x^{6}}\right)_0\right)=A_1$,
$\pr_{\so(n)}\left(R\left(\frac{\partial}{\partial x^{4}},\frac{\partial}{\partial x^{6}}\right)_0\right)=A_2$,
$\pr_{\so(n)}\left(R\left(\frac{\partial}{\partial x^{5}},\frac{\partial}{\partial x^{6}}\right)_0\right)=A_3$
and $\pr_{\so(n)}\left(R\left(\frac{\partial}{\partial x^{1}},\frac{\partial}{\partial x^{6}}\right)_0\right)=
\pr_{\so(n)}\left(R\left(\frac{\partial}{\partial x^{2}},\frac{\partial}{\partial x^{6}}\right)_0\right)=0$.

Note the following.
Let $P\in\Hom(\Real^n,\h)$ be a linear map defined as follows $P(e_1)=P(e_2)=0$, $P(e_3)=A_1$, $P(e_4)=A_2$ and  $P(e_5)=A_3$.
We have $P\in\P(\h)$, $P(\Real^n)=\h$ and
$\pr_{\so(n)}\left(R\left(\frac{\partial}{\partial x^{i}},\frac{\partial}{\partial x^{6}}\right)_0\right)=P(e_i)$ for all $1\leq i\leq 5$.

\section{Main results}

In this section we will construct metrics that for any Riemannian holonomy algebra  $\h\subset\so(n)$
realize the Lie algebras $\g^{1,\h}$, $\g^{2,\h}$, $\g^{3,\h,\varphi}$ and $\g^{4,\h,m,\psi}$
(if  $\g^{3,\h,\varphi}$ and $\g^{4,\h,m,\psi}$ exist)
as holonomy algebras.

Recall that the Lie algebra $\g^{3,\h,\varphi}$ exists only for $\h\subset\so(n)$ with $\z(\h)\neq\{0\}$ and
the Lie algebra $\g^{4,\h,m,\psi}$ exists only for $\h\subset\so(m)$ with $\dim\z(\h)\geq n-m.$

{\bf Constructions of the metrics.}
Let $\h\subset\so(n)$ be the holonomy algebra of a Riemannian manifold.
The Borel-Lichnerowicz theorem (\cite{Bo-Li}) states that we have an orthogonal
decomposition \begin{equation}\label{LM0A}\Real^{n}=\Real^{n_1}\oplus\cdots\oplus\Real^{n_s}\oplus\Real^{n_{s+1}}\end{equation} and the corresponding decomposition into the direct sum of ideals
\begin{equation}\label{LM0B}\h=\h_1\oplus\cdots\oplus\h_s\oplus\{0\}\end{equation} such that $\h$ annulates $\Real^{n_{s+1}}$,
$\h_i(\Real^{n_j})=0$ for $i\neq j$, and $\h_i\subset\so(n_i)$ is an
irreducible subalgebra for $1\leq i\leq s$. Moreover, the Lie algebras $\h_i$ are the holonomy algebras of Riemannian manifolds.
Note that we have (\cite{Le1,G1}) \begin{equation}\label{LM0C}\P(\h)=\P(\h_1)\oplus\cdots\oplus\P(\h_s).\end{equation}
We will assume that the basis $e_1,...,e_n$ of $\Real^n$ is compatible  with the above decomposition of $\Real^n$.

Let $n_0=n_1+\cdots+n_s=n-n_{s+1}$. We see that $\h\subset\so(n_0)$ and $\h$ does not annulate any proper subspace of $\Real^{n_0}$.
Note that in the case of the Lie algebra $\g^{4,\h,m,\psi}$ we have $0<n_0\leq m$.

We will always assume that the indices $b,c,d,f$ run from $0$ to $n+1$,
the indices $i,j,k,l$ run from $1$ to $n$, the indices $\i,\j,\k,\l$ run from $1$ to $n_0$,
the indices $\ii,\jj,\kk,\ll$ run from $n_0+1$ to $n$,  and the indices $\alpha,\beta,\gamma$
run from $1$ to $N$. In case of the Lie algebra  $\g^{4,\h,m,\psi}$ we will also  assume  that
the indices $\iv,\jv,\kv,\lv$ run from $n_0+1$ to $m$ and the indices $\iiv,\jjv,\kkv,\llv$ run from $m+1$ to $n$.
We will use the Einstein rule for sums.

Let $(P_\alpha)_{\alpha=1}^N$ be linearly independent elements of $\P(\h)$ such that the subset
$\{P_\alpha(u)|1\leq\alpha\leq N,\,u\in \Real^{n_0}\}\subset\h$ generates the Lie algebra $\h$. For example, it can be any basis
of the vector space $\P(\h)$.
We have $P_\alpha|_{R^{n_{s+1}}}=0$ and $P_\alpha$ can be considered as linear maps $P_\alpha:\Real^{n_0}\to\h\subset\so(n_0)$.
For each $P_\alpha$ define the numbers $P_{\alpha\j\i}^\k$ such that $P_\alpha(e_\i)e_\j=P_{\alpha\j\i}^\k e_\k.$
Since $P_\alpha\in\P(\h)$, we have
\begin{equation}\label{LM1} P_{\alpha\k\i}^\j=-P_{\alpha\j\i}^\k \text{ and }
 P_{\alpha\j\i}^\k+P_{\alpha\k\j}^\i+P_{\alpha\i\k}^\j=0.\end{equation}

Define the following numbers
\begin{equation}\label{LM2A} a_{\alpha\j\i}^\k=\frac{1}{3\cdot(\alpha-1)!}\left(P_{\alpha\j\i}^\k+P_{\alpha\i\j}^\k\right).\end{equation}
We have \begin{equation}\label{LM3} a_{\alpha\j\i}^\k=a_{\alpha\i\j}^\k.\end{equation}
From \eqref{LM1} it follows that
\begin{equation}\label{LM4} P_{\alpha\j\i}^\k=(\alpha-1)!\left(a_{\alpha\j\i}^\k-a_{\alpha\k\i}^\j\right) \text{ and }
 a_{\alpha\j\i}^\k+a_{\alpha\k\j}^\i+a_{\alpha\i\k}^\j=0.\end{equation}

Let $x^0,...,x^{n+1}$ be the standard coordinates on $\Real^{n+2}$. Consider the following metric
\begin{equation}\label{LM6} g=2dx^0dx^{n+1}+\sum^{n}_{i=1}(dx^i)^2+2\sum^{n_0}_{\i=1}u^\i dx^\i dx^{n+1}+f\cdot(dx^{n+1})^2,\end{equation}
where \begin{equation}\label{LM6A}u^\i=a_{\alpha\j\k}^\i x^\j x^\k(x^{n+1})^{\alpha-1}\end{equation}
and $f$ is a function that will depend on the type of the holonomy algebra that we wish to obtain.

For the Lie algebra $\g^{3,\h,\varphi}$ (if it exists)  define the numbers
\begin{equation}\label{LM7}\varphi_{\alpha\i}=\frac{1}{(\alpha-1)!}\varphi(P_{\alpha}(e_\i)).\end{equation}

For the Lie algebra $\g^{4,\h,m,\psi}$ (if it exists)  define the numbers $\psi_{\alpha\i\iiv}$ such that

\begin{equation}\label{LM8}\frac{1}{(\alpha-1)!}\psi(P_{\alpha}(e_\i))=\sum^{n}_{\iiv=m+1}\psi_{\alpha\i\iiv}e_\iiv.\end{equation}

Suppose that $f(0)=0$, then $g_0=\eta$ and we can identify the tangent space to $\Real^{n+2}$ at $0$ with
the vector space $\Real^{1,n+1}$.

\begin{theorem}\label{LMth1}
The holonomy algebra $\hol_0$ of the metric $g$ at the point $0\in \Real^{n+2}$
depends on the function $f$  in the following way
$$
\begin{array}{ccl}\hline f& &\hol_0 \\\hline
(x^0)^2+\sum_{\ii=n_0+1}^{n}(x^\ii)^2& &\g^{1,\h}\\
\sum_{\ii=n_0+1}^{n}(x^\ii)^2& &\g^{2,\h}\\
2x^0\varphi_{\alpha\i}x^\i(x^{n+1})^{\alpha-1}+\sum_{\ii=n_0+1}^{n}(x^\ii)^2& &\g^{3,\h,\varphi} (\text{if } \z(\h)\neq\{0\}) \\
2\psi_{\alpha\i\iiv}x^\i x^\iiv(x^{n+1})^{\alpha-1}+\sum_{\iv=n_0+1}^{m}(x^\iv)^2& & \g^{4,\h,m,\psi}(\text{if } \dim\z(\h)\geq n-m)\\\hline\end{array}$$

\end{theorem}

As the corollary we get the classification of the weakly-irreducible not irreducible Lorentzian holonomy algebras.

\begin{theorem}\label{LMth2}
A weakly-irreducible not irreducible subalgebra  $\g\subset\so(1,n+1)$ is the holonomy algebra of a Lorentzian manifold if and only if $\g$ is conjugated to one
of the subalgebras $\g^{1,\h},\g^{2,\h},\g^{3,\h,\varphi},\g^{4,\h,m,\psi}\subset\so(1,n+1)_{\Real p}$, where $\h\subset\so(n)$ is the
holonomy algebra of a Riemannian manifold.
\end{theorem}

From theorem \ref{LMth2}, Wu's theorem and Berger's list it follows that
{\it the holonomy algebra $\hol\subset\so(1,N+1)$ of any Lorentzian manifold of dimension $N+2$
has the form $\hol=\g\oplus\h_1\oplus\cdots\oplus\h_r$,
where either $\g=\so(1,n+1)$ or $\g$ is a Lie algebra from theorem \ref{LMth2}, and $\h_i\subset\so(n_i)$ are the irreducible holonomy
algebras of Riemannian manifolds ($N=n+n_1+\cdots+n_r$).}

{\bf Explanation of the idea of the constructions.}
Now we compare our method of constructions with the example of A. Ikemakhen.
Let us construct the metric for $\g^{2,\rho(\so(3))}\subset\so(1,6)$ by
our method. Take $P\in\P(\rho(\so(3)))$ defined as $P(e_1)=P(e_2)=0$, $P(e_3)=A_1$, $P(e_4)=A_2$ and $P(e_5)=A_3$.
By our constructions, we have
$$g=2dx^0dx^{6}+\sum^{5}_{i=1}(dx^i)^2+2\sum^{5}_{i=1}u^i dx^{i} dx^{6},$$
where $$\begin{array}{ll}u^1=-\frac{2}{3}((x^3)^2+4(x^4)^2+(x^5)^2),& u^2=\frac{2\sqrt{3}}{3}((x^3)^2-(x^5)^2),\\
u^3=\frac{2}{3}(x^1x^3-\sqrt{3}x^2x^3-3x^4x^5-(x^5)^2),& u^4=\frac{8}{3}x^1x^4,\\
u^5=\frac{2}{3}(x^1x^5+\sqrt{3}x^2x^5+3x^3x^4+x^3x^5).& \end{array}$$
We still have $\pr_{\so(n)}\left(R\left(\frac{\partial}{\partial x^{3}},\frac{\partial}{\partial x^{6}}\right)_0\right)=A_1$,
$\pr_{\so(n)}\left(R\left(\frac{\partial}{\partial x^{4}},\frac{\partial}{\partial x^{6}}\right)_0\right)=A_2$,
$\pr_{\so(n)}\left(R\left(\frac{\partial}{\partial x^{5}},\frac{\partial}{\partial x^{6}}\right)_0\right)=A_3$
and $\pr_{\so(n)}\left(R\left(\frac{\partial}{\partial x^{1}},\frac{\partial}{\partial x^{6}}\right)_0\right)=
\pr_{\so(n)}\left(R\left(\frac{\partial}{\partial x^{2}},\frac{\partial}{\partial x^{6}}\right)_0\right)=0$.

The reason why we obtain another metric is the following. The idea of our constructions is to find the constants $a_{\alpha\j\i}^\k$ such that
$$\begin{array}{l}
\pr_{\so(n)}\left(R\left(\frac{\partial}{\partial x^{\i}},\frac{\partial}{\partial x^{n+1}}\right)_0\right)=P_1(e_\i),\\
...,\\
\pr_{\so(n)}\left(\nabla^{N-1}  R\left(\frac{\partial}{\partial x^{\i}},\frac{\partial}{\partial x^{n+1}};
\frac{\partial}{\partial x^{n+1}};\cdots ;\frac{\partial}{\partial x^{n+1}}\right)_0\right)=P_N(e_\i).
\end{array}$$
These conditions give us the system of equations
$$\left\{\begin{array}{rl}(\alpha-1)!\left(a_{\alpha\j\i}^\k-a_{\alpha\k\i}^\j\right)=&P_{\alpha\j\i}^\k,\\
a_{\alpha\j\i}^\k-a_{\alpha\i\j}^\k=&0.\end{array}\right.$$
One of the solutions of this system is given by \eqref{LM2A}, but this system can have other solutions.
In the example we use the solution given by \eqref{LM2A},
taking another solution of the above system, we can obtain the metric constructed by
A. Ikemakhen.

Thus the choose of the functions $\u^\i$ given by \eqref{LM6A} guarantees us that the orthogonal part of the holonomy algebra
$\hol_0$ coincides with the given Riemannian  holonomy algebra $\h\subset\so(n)$ (the other values of
$\pr_{\so(n)}(\nabla^{r}  R)$ does not give us any thing new). This also guarantees us the inclusion $\Real^{n_0}\subset\hol_0$.
The reason why we choose the function $f$ as in theorem \ref{LMth1} can be easily understood from
the following formulas

$$\pr_{\Real}\left(R\left(\frac{\partial}{\partial x^{0}},\frac{\partial}{\partial x^{n+1}}\right)_0\right)=
\frac{1}{2}\frac{\partial^2 f}{(\partial x^{0})^2} \text{ (we use this for } \g^{1,\h}),$$
$$\pr_{\Real}\left(\nabla^{\alpha-1}  R\left(\frac{\partial}{\partial x^{\i}},\frac{\partial}{\partial x^{n+1}};
\frac{\partial}{\partial x^{n+1}};\cdots ;\frac{\partial}{\partial x^{n+1}}\right)_0\right)=
\frac{1}{2}\frac{\partial^{\alpha+1} f}{\partial x^{0}\partial x^{\i}(\partial x^{n+1})^{\alpha-1}} \text{ (we use this for }
\g^{3,\h,\varphi}),$$
$$\pr_{\Real^n}\left(\nabla^{\alpha-1}  R\left(\frac{\partial}{\partial x^{\ii}},\frac{\partial}{\partial x^{n+1}};
\frac{\partial}{\partial x^{n+1}};\cdots ;\frac{\partial}{\partial x^{n+1}}\right)_0\right)=
\frac{1}{2}\sum_{\jj=n_0+1}^{n}
\frac{\partial^{\alpha+1} f}{\partial x^{\ii}\partial x^{\jj}(\partial x^{n+1})^{\alpha-1}}e_\jj
$$
$$\text{ (for $\alpha=0$ we use this for all algebras, for $\alpha\geq 0$ we use this for }
\g^{4,\h,m,\psi}).$$

Es application we construct  metrics for the Lie algebras $\g^{2,g_2}\subset\so(1,8)$ and $\g^{2,\spin(7)}\subset\so(1,9)$.

{\bf Example (Metric with the holonomy algebra $\g^{2,g_2}\subset\so(1,8)$).}
Consider the Lie subalgebra $g_2\subset\so(7)$. The vector subspace $g_2\subset\so(7)$ is spanned by the following matrices
(\cite{HelgaInes})
$$\begin{array}{lllll}
A_1=E_{12}-E_{34},&A_2=E_{12}-E_{56},&A_3=E_{13}+E_{24},&A_4=E_{13}-E_{67},&A_5=E_{14}-E_{23},\\
A_6=E_{14}-E_{57},&A_7=E_{15}+E_{26},&A_8=E_{15}+E_{47},&A_9=E_{16}-E_{25},&A_{10}=E_{16}+E_{37},\\
A_{11}=E_{17}-E_{36},&A_{12}=E_{17}-E_{45},&A_{13}=E_{27}-E_{35},&A_{14}=E_{27}+E_{46},&
\end{array}$$ where $E_{ij}\in\so(7)$ ($i<j$) is the
skew-symmetric matrix such that $(E_{ij})_{ij}=1$, $(E_{ij})_{ji}=-1$ and $(E_{ij})_{kl}=0$ for other $k$ and $l$.

Consider the linear map $P\in\Hom(\Real^7,g_2)$ defined as
$$\begin{array}{llll}P(e_1)=A_6,&P(e_2)=A_4+A_5,&P(e_3)=A_1+A_7,&P(e_4)=A_1,\\P(e_5)=A_4,&P(e_6)=-A_5+A_6,&P(e_7)=A_7.&\end{array}$$
It can be checked that $P\in\P(g_2)$. Moreover,  the elements $A_1,A_4,A_5,A_6,A_7\in g_2$ generate the Lie algebra
$g_2$.

The holonomy algebra of the metric
$$g=2dx^0dx^{8}+\sum^{7}_{i=1}(dx^i)^2+2\sum^{7}_{i=1}u^i dx^{i} dx^{8},$$
where
$$\begin{array}{ll}
u^1=\frac{2}{3}(2x^2x^3+x^1x^4+2x^2x^4+2x^3x^5+x^5x^7),&\\u^2=\frac{2}{3}(-x^1x^3-x^2x^3-x^1x^4+2x^3x^6+x^6x^7),&\\
u^3=\frac{2}{3}(-x^1x^2+(x^2)^2-x^3x^4-(x^4)^2-x^1x^5-x^2x^6),&u^4=\frac{2}{3}(-(x^1)^2-x^1x^2+(x^3)^2+x^3x^4),\\
u^5=\frac{2}{3}(-x^1x^3-2x^1x^7-x^6x^7),&u^6=\frac{2}{3}(-x^2x^3-2x^2x^7-x^5x^7),\\
u^7=\frac{2}{3}(x^1x^5+x^2x^6+2x^5x^6),&\end{array}$$
at the point $0\in\Real^9$ is $\g^{2,g_2}\subset\so(1,8)$.

Using computer it can be checked that $\dim\P(g_2)=64$. This means that we can construct quite a big number of
metrics with the holonomy algebra $\g^{2,g_2}\subset\so(1,8)$.

{\bf Example (Metric with the holonomy algebra $\g^{2,\spin(7)}\subset\so(1,9)$).}
Consider the Lie subalgebra $\spin(7)\subset\so(8)$.
The vector subspace $\spin(7)\subset\so(8)$ is spanned by the following matrices
(\cite{HelgaInes})
$$\begin{array}{lllll}
A_{1}=E_{12}+E_{34},&A_{2}=E_{13}-E_{24},&A_{3}=E_{14}+E_{23},&A_{4}=E_{56}+E_{78},&A_{5}=-E_{57}+E_{68},\\
A_{6}=E_{58}+E_{67},&A_{7}=-E_{15}+E_{26},&A_{8}=E_{12}+E_{56},&A_{9}=E_{16}+E_{25},&A_{10}=E_{37}-E_{48},\\
A_{11}=E_{38}+E_{47},&A_{12}=E_{17}+E_{28},&A_{13}=E_{18}-E_{27},&A_{14}=E_{35}+E_{46},&A_{15}=E_{36}-E_{45},\\
A_{16}=E_{18}+E_{36},&A_{17}=E_{17}+E_{35},&A_{18}=E_{26}-E_{48},&A_{19}=E_{25}+E_{38},&A_{20}=E_{23}+E_{67},\\
A_{21}=E_{24}+E_{57}.
\end{array}$$
Consider the linear map $P\in\Hom(\Real^8,\spin(7))$ defined as
$$\begin{array}{llll}
P(e_1)=0,&P(e_2)=-A_{14},&P(e_3)=0,&P(e_4)=A_{21},\\P(e_5)=A_{20},&P(e_6)=A_{21}-A_{18},&P(e_7)=A_{15}-A_{16},&P(e_7)=A_{14}-A_{17}.
\end{array}$$
It can be checked that $P\in\P(\spin(7))$. Moreover,  the elements $A_{14},A_{15}-A_{16},A_{17},A_{18},A_{20},A_{21}\in  \spin(7)$
generate the Lie algebra
$\spin(7)$.

The holonomy algebra of the metric
$$g=2dx^0dx^{9}+\sum^{8}_{i=1}(dx^i)^2+2\sum^{8}_{i=1}u^i dx^{i} dx^{9},$$
where
$$\begin{array}{ll}
u^1=-\frac{4}{3}x^7x^8,&
u^2=\frac{2}{3}((x^4)^2+x^3x^5+x^4x^6-(x^6)^2),\\
u^3=-\frac{4}{3}x^2x^5,&
u^4=\frac{2}{3}(-x^2x^4-2x^2x^6-x^5x^7+2x^6x^8),\\
u^5=\frac{2}{3}(x^2x^3+2x^4x^7+x^6x^7),&
u^6=\frac{2}{3}(x^2x^4+x^2x^6+x^5x^7-x^4x^8),\\
u^7=\frac{2}{3}(-x^4x^5-2x^5x^6+x^1x^8),&
u^8=\frac{2}{3}(-x^4x^6+x^1x^7),
\end{array}$$
at the point $0\in\Real^9$ is $\g^{2,\spin(7)}\subset\so(1,9)$.

Note that $\dim\P(\spin(7))=112$.

{\bf Proof of theorem \ref{LMth1}.}
We will prove the theorem for the case of algebras of type 4. For other types the proof is analogous.

Since the coefficients of the metric $g$ are polynomial functions, the Levi-Civita
connection given by $g$ is analytic and the Lie algebra $\hol_0$ is generated by the operators
$$R(X,Y)_0,\nabla R(X,Y;Z_1)_0,\nabla^2 R(X,Y;Z_1;Z_2)_0,...\in\so(1,n+1),$$
where  $\nabla^r R(X,Y;Z_1;...;Z_r)=(\nabla_{Z_r}\cdots\nabla_{Z_1}R)(X,Y)$ and
$X$, $Y$, $Z_1$, $Z_2$,... are vectors at the point $0$.

To prove the theorem we will compute the components of the curvature tensor and of its derivatives.

The non-zero  Christoffel symbols for the metric $g$ given by \eqref{LM6} are the following
$$\begin{array}{rl}
\Ga^0_{0n+1}&=\frac{1}{2}\frac{\partial f}{\partial x^0},\\
\Ga^0_{n+1\, n+1}&=\frac{1}{2}\left(-\sum_{\i=1}^{n_0}u^\i\left(2\frac{\partial u^\i}{\partial x^{n+1}}-\frac{\partial f}{\partial x^{\i}}\right)+
\frac{\partial f}{\partial x^{n+1}}\right)+\frac{\partial f}{\partial x^{0}}\left(f-\sum_{\i=1}^{n_0}(u^\i)^2\right),\\
\Ga^0_{\i\j}&=\frac{1}{2}\left(\frac{\partial u^\i}{\partial x^\j}+\frac{\partial u^\j}{\partial x^\i}\right),\\
\Ga^0_{\i n+1}&=\frac{1}{2}\sum_{\j=1}^{n_0}u^\j\left(\frac{\partial u^\i}{\partial x^\j}-\frac{\partial u^\j}{\partial x^\i}\right)+
\frac{1}{2}\frac{\partial f}{\partial x^\i},\\
\Ga^0_{\ii n+1}&=\frac{1}{2}\frac{\partial f}{\partial x^\ii},\\
\Ga^\i_{\j n+1}&=\frac{1}{2}\left(\frac{\partial u^\i}{\partial x^\j}-\frac{\partial u^\j}{\partial x^\i}\right),\\
\Ga^\i_{n+1\, n+1}&=\frac{1}{2}\left(-2\frac{\partial u^\i}{\partial x^{n+1}}-\frac{\partial f}{\partial x^\i}+u^\i\frac{\partial f}{\partial x^0}\right),\\
\Ga^\ii_{n+1\, n+1}&=-\frac{1}{2}\frac{\partial f}{\partial x^\ii},\\
\Ga^{n+1}_{n+1\, n+1}&=-\frac{1}{2}\frac{\partial f}{\partial x^0}.\end{array}$$

Suppose that $\dim\z(\h)\geq n-m$.
Let $f=2\psi_{\alpha\i\iiv}x^\i x^\iiv(x^{n+1})^{\alpha-1}+\sum_{\iv=n_0+1}^{m}(x^\iv)^2$. We must prove that $\hol_0= \g^{4,\h,m,\psi}$.
We have the following non-zero Christoffel symbols

\begin{multline}
\Ga^0_{n+1\, n+1}=-\sum_{\i=1}^{n_0}u^\i((\alpha-1)a^\i_{\alpha\j\k}x^\j x^\k(x^{n+1})^{\alpha-2}-
\psi_{\alpha\i\iiv}x^\iiv(x^{n+1})^{\alpha-1})+\\(\alpha-1)\psi_{\alpha\i\iiv}x^\i x^\iiv(x^{n+1})^{\alpha-2}, \label{LM8E}
\end{multline}
\begin{align}
\Ga^0_{\i\j}&=\frac{1}{2}(a^\i_{\alpha\j\k}+a^\j_{\alpha\i\k})x^\k(x^{n+1})^{\alpha-1}, \label{LM8A}\\
\Ga^0_{\i,n+1}&=\frac{1}{(\alpha-1)!}u^\j P^\i_{\alpha\j\k}x^\k(x^{n+1})^{\alpha-1}+
\psi_{\alpha\i\iiv}x^\iiv(x^{n+1})^{\alpha-1}, \label{LM8B}\\
\Ga^\i_{\j n+1}&=\frac{1}{(\alpha-1)!}P^\i_{\alpha\j\k}x^\k(x^{n+1})^{\alpha-1}, \label{LM8F}\\
\Ga^0_{\iv n+1}&=x^\iv, \label{LM8C}\\
\Ga^0_{\iiv n+1}&=\psi_{\alpha\i\iiv}x^\i(x^{n+1})^{\alpha-1}, \label{LM8D}\\
\Ga^\iv_{n+1\, n+1}&=-x^\iv,\label{LM8D1} \\
\Ga^\iiv_{n+1\, n+1}&=-\psi_{\alpha\i\iiv}x^\i(x^{n+1})^{\alpha-1}.\label{LM8D2}
\end{align}

In particular, note the following
\begin{equation}\Ga^k_{ij}=\Ga^{n+1}_{bc}=\Ga^b_{0c}=\Ga^i_{\ii n+1}=0. \label{LM2226}
\end{equation}

For $r\geq 0$ let $R^b_{c,d,f;f_1;...;f_r}$ be the functions such that
$\nabla^r R(\frac{\partial}{\partial x^d},\frac{\partial}{\partial x^f};
\frac{\partial}{\partial x^{f_1}};\cdots; \frac{\partial}{\partial x^{f_r}})\frac{\partial}{\partial x^c}=
R^b_{c,d,f;f_1;...;f_r}\frac{\partial}{\partial x^b}$.

One can compute the following components of the curvature tensor
\begin{align}
R^\k_{\j\i\,n+1}&=\frac{1}{(\alpha-1)!}P^\k_{\alpha\j\i}(x^{n+1})^{\alpha-1}, R^\kk_{j\i\,n+1}=0, \label{LM9A}\\
R^k_{jbc}&=0 \text{ if } (b,c)\notin\{1,...,n_0\}\times\{n+1\}\cup\{n+1\}\times\{1,...,n_0\},      \label{LM9A1}\\
R^0_{\iiv\i n+1}&=\psi_{\alpha\i\iiv}(x^{n+1})^{\alpha-1},\label{LM9C}\\
R^0_{\iiv bc}&=0\text{ if } (b,c)\notin\{1,...,n_0\}\times\{n+1\}\cup\{n+1\}\times\{1,...,n_0\},      \label{LM9C1}\\
R^0_{\k\i\j}&=\frac{1}{(\alpha-1)!}P^\j_{\alpha\i\k}(x^{n+1})^{\alpha-1}, \label{LM9D}\\
R^0_{\iv\iv n+1}&=1,R^\iv_{n+1\,\iv n+1}=-1,R^b_{c\iv n+1}=0 \text{ if } (b,c)\neq(0,\iv)\text{ and } (b,c)\neq(\iv,n+1),\label{LM9E}\\
R^0_{0bc}&=0. \label{LM9E1}
\end{align}

Using \eqref{LM2226}, we get
\begin{equation}\label{LM10}
\Ga^h_{b f}R^f_{cdg}=\Ga^h_{b i}R^i_{cdg},\quad \Ga^f_{b h}R^c_{fdg}=\Ga^i_{b h}R^c_{i dg},\quad \Ga^f_{b h}R^c_{dfg}=\Ga^i_{b h}R^c_{dig},
\quad \Ga^f_{b h}R^c_{dgf}=\Ga^i_{b h}R^c_{ dgi}.\end{equation}

From equalities \eqref{LM8E} -- \eqref{LM8D2} it follows that
\begin{equation}\label{LM12}
\Ga^i_{bc}R^d_{ief}=\Ga^\i_{bc}R^d_{\i ef},\quad \Ga^i_{bc}R^d_{eif}=\Ga^\i_{bc}R^d_{e\i f}, \quad \Ga^i_{bc}R^d_{efi}=\Ga^\i_{bc}R^d_{ef\i}
\text{ if } b\neq n+1 \text{ or } c\neq n+1,\end{equation}
\begin{equation}\label{LM13}\Ga^b_{ic}R^i_{def}=\Ga^b_{\i c}R^\i_{def}\text{ if } b\neq 0 \text{ or } c\neq n+1.\end{equation}

{\bf Proof of the inclusion $\g^{4,\h,m,\psi}\subset\hol_0$.}


\begin{lem}\label{LMlem1}
For any $1\leq r\leq N$ we have

\begin{itemize}

\item[{\rm 1)}] $R^\k_{\j\i n+1;\underbrace{n+1;\cdots;n+1}_{r-1 \text{ times}}}=
\sum_{\alpha=r}^N\frac{1}{(\alpha-r)!}P^\k_{\alpha\j\i}(x^{n+1})^{\alpha-r}+y^\k_{\alpha\j\i r}$, where
$y^\k_{\alpha\j\i r}$ are functions such that $y^\k_{\alpha\j\i r}(0)=\left(\frac{\partial y^\k_{\alpha\j\i r}}{\partial x^{n+1}}\right)(0)=\cdots=0$;

\item[{\rm 2)}] $R^0_{\iiv\i n+1;\underbrace{n+1;\cdots;n+1}_{r-1 \text{ times}}}=\sum_{\alpha=r}^N\frac{(\alpha-1)!}{(\alpha-r)!}\psi_{\alpha\i\iiv}
(x^{n+1})^{\alpha-r}+z^\k_{\alpha\j\i r}$, where
$z^\k_{\alpha\j\i r}$ are functions such that $z^\k_{\alpha\j\i r}(0)=\left(\frac{\partial z^\k_{\alpha\j\i r}}{\partial x^{n+1}}\right)(0)=\cdots=0$.
\end{itemize}
\end{lem}

{\it Proof.} We will prove this lemma by induction over $r$.
For $r=1$ the lemma follows from \eqref{LM9A} and \eqref{LM9C}. Fix $r_0>1$  and assume that the lemma is true for all $r<r_0$.
We must prove that the lemma is true for $r=r_0$. We have \\
$R^\k_{\j\i n+1;\underbrace{n+1;\cdots;n+1}_{r \text{ times}}}=$
$$\begin{array}{rl}
&\frac{\partial R^\k_{\j\i n+1;n+1;\cdots;n+1(r-1 \text{ times})}}{\partial x^{n+1}}
+\Ga^\k_{n+1\, l}R^l_{\j\i n+1;\underbrace{n+1;\cdots;n+1}_{r-1 \text{ times}}}\\&
-\Ga^l_{n+1\, \j}R^\k_{l\i n+1;\underbrace{n+1;\cdots;n+1}_{r-1 \text{ times}}}
-\Ga^l_{n+1\, \i}R^\k_{\j l n+1;\underbrace{n+1;\cdots;n+1}_{r-1 \text{ times}}}\\&
-\Ga^l_{n+1\, n+1}R^\k_{\j\i l;\underbrace{n+1;\cdots;n+1}_{r-1 \text{ times}}}\\&
-\Ga^l_{n+1\, n+1}R^\k_{\j\i n+1;l;\underbrace{n+1;\cdots;n+1}_{r-2 \text{ times}}}-\cdots
-\Ga^l_{n+1\, n+1}R^\k_{\j\i n+1;\underbrace{n+1;\cdots;n+1}_{r-2 \text{ times}};l}\\
=& \sum_{\alpha=r+1}^N\frac{1}{(\alpha-r-1)!}P^\k_{\alpha\j\i}(x^{n+1})^{\alpha-r-1}+\frac{\partial y^\k_{\alpha\j\i r}}{\partial x^{n+1}}+\tilde{y}^\k_{\alpha\j\i r}\\
=&\sum_{\alpha=r+1}^N\frac{1}{(\alpha-r-1)!}P^\k_{\alpha\j\i}(x^{n+1})^{\alpha-r-1}+y^\k_{\alpha\j\i r+1}.
\end{array}$$
Claim 1) follows from the fact that  all Christoffel symbols and all their derivatives with respect to $x^{n+1}$ vanish at the point $0$.
The proof of claim 2) is analogous. The lemma is proved. $\Box$

From lemma \ref{LMlem1} it follows that for any $1\leq r\leq N$ we have
\begin{align}R^\k_{\j\i n+1;\underbrace{n+1;\cdots;n+1}_{r-1 \text{ times}}}(0)&=P^\k_{r\j\i},\label{LM14}\\
R^0_{\iiv\i n+1;\underbrace{n+1;\cdots;n+1}_{r-1 \text{ times}}}(0)&=(r-1)!\psi_{r\i\iiv}.\label{LM15}\end{align}
Similarly we can prove that from \eqref{LM9A} it follows that
\begin{equation}\label{LM16}R^\kk_{i\i n+1;\underbrace{n+1;\cdots;n+1}_{r-1 \text{ times}}}(0)=0.\end{equation}
From \eqref{LM8},  \eqref{LM14},  \eqref{LM15} and   \eqref{LM16} it follows that
\begin{equation}\label{LM17}R({e_\i, e_{n+1};\underbrace{e_{n+1};\cdots;e_{n+1}}_{r-1 \text{ times}}})_0=(0,P_r(e_\i),X_{r\i}+\psi(P_r(e_\i))),\end{equation}
where $X_{r\i}\in\spa\{e_1,...,e_m\}.$
Since the elements $P_r(e_\i)$ generate the Lie algebra $\h$ and $\pr_{\so(n)}\hol_0$ is a Lie algebra,  from \eqref{LM17} we get
\begin{equation}\label{LM18}\h\subset\pr_{\so(n)}\hol_0.\end{equation}
Using \eqref{LM9D} we can prove that for all $1\leq r\leq N$ we have
$$R^0_{\k\i\j;\underbrace{n+1;\cdots;n+1}_{r-1 \text{ times}}}(0)=P^\j_{r\i\k}.$$
This means that
\begin{equation}\label{LM19} (0,0,P^\j_{r\i\k}e_\j)\in\hol_0.\end{equation}
Recall that we have decompositions \eqref{LM0A}, \eqref{LM0B} and \eqref{LM0C}.
Since $\h$ is generated by the images of the elements $P_\alpha$,
for any $2\leq i\leq s$ there exist $\alpha,\i,\j,\k$ such that $n_1+\cdots+n_{i-1}+1\leq\j\leq n_1+\cdots+n_i$
and $P^\j_{\alpha\i\k}\neq 0$.
Combining this with \eqref{LM19}, we get
$$\{(0,0,X)|X\in\Real^{n_i}\}\cap \hol_0\neq\{0\}.$$
Since $\h_i\subset\so(n_i)$ is an irreducible subalgebra, $\h_i\subset\pr_{\so(n)}\hol_0$, and for any $A\in\h$, $Z\in\Real^n$, $Y\in\Real^{n_i}$ holds
$$[(0,A,Z),(0,0,Y)]=(0,0,AY)\in\{(0,0,X)|X\in\Real^{n_i}\}\cap \hol_0,$$
we see that
$$\{(0,0,X)|X\in\Real^{n_i}\}\subset \hol_0,$$
hence,
\begin{equation}\label{LM20}\{(0,0,X)|X\in\spa\{e_1,...,e_{n_0}\}\}\subset \hol_0.\end{equation}
From \eqref{LM9E} it follows that $R(e_\iv,q)=(0,0,e_\iv)$, hence,  \begin{equation}\label{LM21}\{(0,0,X)|X\in\spa\{e_{n_0+1},...,e_{m}\}\}\subset \hol_0.\end{equation}
From \eqref{LM17}, \eqref{LM20}, \eqref{LM21} and the fact that $\h$ is generated by the elements $P_\alpha(e_\i)$ it follows that
$$\g^{4,\h,m,\psi}\subset\hol_0.$$

{\bf Proof of the inclusion $\hol_0\subset\g^{4,\h,m,\psi}$.}

\begin{lem}\label{LMlem2}

For any $r\geq 0$ we have
\begin{itemize}

\item[{\rm 1)}] $R^k_{jil;f_1;\cdots;f_r}=0$;

\item[{\rm 2)}] $R^\ii_{jin+1;f_1;\cdots;f_r}=0$;

\item[{\rm 3)}] $R^0_{\iiv il;f_1;\cdots;f_r}=0$;

\item[{\rm 4)}] $R^\k_{\j \ii n+1;f_1;\cdots;f_r}=0$;

\item[{\rm 5)}] $R^0_{\iiv\ii n+1;f_1;\cdots;f_r}=0$;

\item[{\rm 6)}] $R^0_{0bc;f_1;\cdots;f_r}=0$;

\item[{\rm 7)}] $R^\k_{\j\i n+1;f_1;\cdots;f_r}=\sum_{t\in T_{\i f_1...f_r}}z_tA^\k_{t\j}$, where $T_{\i f_1...f_r}$ is a finite set of indeces,
$z_t$ are functions and $A_t\in\h$.

\item[{\rm 8)}] $R^0_{\iiv\i n+1;f_1;\cdots;f_r}=\sum_{t\in T_{\i f_1...f_r}}z_t\psi _{t\iiv}$, where $\psi _{t\iiv}$ are numbers such that
$\psi(A_t)=\sum_{\iiv=m+1}^{n} \psi _{t\iiv}e_\iiv$.
\end{itemize}\end{lem}

{\it Proof.} We will prove the claims of the lemma by induction over $r$. For $r=0$ the claims follow from \eqref{LM9A} -- \eqref{LM9E1}.
Fix $r>0$  and assume that the lemma is true for  $r$. We must prove that the lemma is true for $r+1$.

1) We have\\
$R^k_{jil;f_1;\cdots;f_r;l_1}=$
$$\begin{array}{rl}
&\frac{\partial R^k_{jil;f_1;\cdots;f_r}}{\partial x^{l_1}}+
\Ga^k_{l_1l_2}R^{l_2}_{jil;f_1;\cdots;f_r} -\Ga^{l_2}_{l_1j}R^{k}_{l_2il;f_1;\cdots;f_r}\\
&-\Ga^{l_2}_{l_1i}R^{k}_{jl_2l;f_1;\cdots;f_r}-\Ga^{l_2}_{l_1l}R^{k}_{jil_2;f_1;\cdots;f_r}
-\Ga^{l_2}_{l_1f_1}R^{k}_{jil;l_2;f_2:\cdots;f_r}-\cdots-\Ga^{l_2}_{l_1f_r}R^{k}_{jil;f_1:\cdots;f_{r-1};l_2}.
\end{array}$$
From \eqref{LM2226} and the inductive hypothesis it follows that
$$R^k_{jil;f_1;\cdots;f_r;l_1}=0.$$

Similarly,\\ $R^k_{jil;f_1;\cdots;f_r;n+1}=$
$$\begin{array}{rl} &\frac{\partial R^k_{jil;f_1;\cdots;f_r}}{\partial x^{n+1}}+
\Ga^k_{n+1\,l_2}R^{l_2}_{jil;f_1;\cdots;f_r} -\Ga^{l_2}_{n+1\,j}R^{k}_{l_2il;f_1;\cdots;f_r}-\Ga^{l_2}_{n+1\,i}R^{k}_{jl_2l;f_1;\cdots;f_r}\\
&-\Ga^{l_2}_{n+1\,l}R^{k}_{jil_2;f_1;\cdots;f_r}
-\Ga^{l_2}_{n+1\,f_1}R^{k}_{jil;l_2;f_2:\cdots;f_r}-\cdots-\Ga^{l_2}_{n+1\,f_r}R^{k}_{jil;f_1:\cdots;f_{r-1};l_2}=0.
\end{array}$$

Claim 1) is proved. The proofs of claims 2) -- 6) are analogous.

7) We have\\ $R^\k _{\j \i n+1;f_1;\cdots;f_r;l}=$
$$\begin{array}{rl} &\frac{\partial R^\k _{\j \i n+1;f_1;\cdots;f_r}}{\partial x^{l}}+
\Ga^\k _{ll_1}R^{l_1}_{\j \i n+1;f_1;\cdots;f_r} -\Ga^{l_1}_{l\j }R^{\k }_{l_1\i n+1;f_1;\cdots;f_r}
-\Ga^{l_1}_{l\i }R^{\k }_{\j l_1n+1;f_1;\cdots;f_r}\\
&-\Ga^{l_1}_{ln+1}R^{\k }_{\j \i l_1;f_1;\cdots;f_r}
-\Ga^{l_1}_{lf_1}R^{\k }_{\j \i n+1;l_1;f_2:\cdots;f_r}-\cdots-\Ga^{l_1}_{lf_r}R^{\k }_{\j \i n+1;f_1:\cdots;f_{r-1};l_1}.
\end{array}$$
From \eqref{LM2226}, claim 1) of the lemma and the inductive hypothesis it follows that
\begin{multline}\label{LM22}
R^\k _{\j \i n+1;f_1;\cdots;f_r;l}=\\ \sum_{t\in T_{\i f_1...f_r}}\frac{\partial z_t}{\partial x^l}A^\k_{t\j}-
\sum_{l_1=1}^{n}\sum_{t\in T_{\i l_1f_2...f_r}}z_t\Ga^{l_1}_{lf_1}A^\k_{t\j}-\cdots-\sum_{l_1=1}^{n}\sum_{t\in T_{\i f_1...f_{r-1}l_1}}z_t\Ga^{l_1}_{lf_r}A^\k_{t\j}.
\end{multline}

Furthermore, $R^\k _{\j \i n+1;f_1;\cdots;f_r;n+1}=$
$$\begin{array}{rl} &\frac{\partial R^\k _{\j \i n+1;f_1;\cdots;f_r}}{\partial x^{n+1}}+
\Ga^\k _{n+1\,l_1}R^{l_1}_{\j \i n+1;f_1;\cdots;f_r} -\Ga^{l_1}_{n+1\,\j }R^{\k }_{l_1\i n+1;f_1;\cdots;f_r}-\Ga^{l_1}_{n+1\,\i }R^{\k }_{\j l_1n+1;f_1;\cdots;f_r}\\
&-\Ga^{l_1}_{n+1\,n+1}R^{\k }_{\j \i l_1;f_1;\cdots;f_r}
-\Ga^{l_1}_{n+1\,f_1}R^{\k }_{\j \i n+1;l_1;f_2:\cdots;f_r}-\cdots-\Ga^{l_1}_{n+1\,f_r}R^{\k }_{\j \i n+1;f_1:\cdots;f_{r-1};l_1}.
\end{array}$$

From \eqref{LM8F} and \eqref{LM2226} it follows that
$$\Ga^\k _{n+1\,l_1}R^{l_1}_{\j \i n+1;f_1;\cdots;f_r} -\Ga^{l_1}_{n+1\,\j }R^{\k }_{l_1\i n+1;f_1;\cdots;f_r}=
\frac{1}{(\alpha-1)!}x^{\l_2}(x^{n+1})^{\alpha-1}[P_\alpha(e_{\l_2}),R(e_\i,e_{n+1};e_{f_1};...;e_{f_r})]^\k_\j.$$

Using this, claim 1) and inductive hypothesis, we get
\begin{multline} \label{LM23}R^\k _{\j \i n+1;f_1;\cdots;f_r;n+1}=\sum_{t\in T_{\i f_1...f_r}}\frac{\partial z_t}{\partial x^{n+1}}A^\k_{t\j}+
\frac{1}{(\alpha-1)!}x^{\l_2}(x^{n+1})^{\alpha-1}[P_\alpha(e_{\l_2}),R(e_\i,e_{n+1};e_{f_1};...;e_{f_r})]^\k_\j\\
-\sum_{l_1=1}^{n}\sum_{t\in T_{ l_1f_1...f_r}}z_t\Ga^{l_1}_{n+1\i}A^\k_{t\j}\\
-\sum_{l_1=1}^{n}\sum_{t\in T_{\i l_1f_2...f_r}}z_t\Ga^{l_1}_{n+1 f_1}A^\k_{t\j}-
\cdots-\sum_{l_1=1}^{n}\sum_{t\in T_{\i f_1...f_{r-1}l_1}}z_t\Ga^{l_1}_{n+1f_r}A^\k_{t\j}.\end{multline}
This proves  claim  7).

8) We have  $$\begin{array}{rl} R^0 _{\iiv \i n+1;f_1;\cdots;f_r;l}=&\frac{\partial R^0 _{\iiv \i n+1;f_1;\cdots;f_r}}{\partial x^{l}}+
\Ga^0 _{ll_1}R^{l_1}_{\iiv \i n+1;f_1;\cdots;f_r} -\Ga^{l_1}_{l\iiv }R^{0 }_{l_1\i n+1;f_1;\cdots;f_r}
-\Ga^{l_1}_{l\i }R^{0 }_{\iiv l_1n+1;f_1;\cdots;f_r}\\
&-\Ga^{l_1}_{ln+1}R^{0 }_{\iiv \i l_1;f_1;\cdots;f_r}
-\Ga^{l_1}_{lf_1}R^{0 }_{\iiv \i n+1;l_1;f_2:\cdots;f_r}-\cdots-\Ga^{l_1}_{lf_r}R^{0 }_{\iiv \i n+1;f_1:\cdots;f_{r-1};l_1}.
\end{array}$$

Using this,  \eqref{LM2226},  claims 2) and 3), and the inductive hypothesis, we get
\begin{equation}\label{LM24} R^0_{\iiv \i n+1;f_1;\cdots;f_r;l}=\sum_{t\in T_{\i f_1...f_r}}\frac{\partial z_t}{\partial x^l}\psi_{t\iiv}-
\sum_{l_1=1}^{n}\sum_{t\in T_{\i l_1f_2...f_r}}z_t\Ga^{l_1}_{lf_1}\psi_{t\iiv}-\cdots-\sum_{l_1=1}^{n}\sum_{t\in T_{\i f_1...f_{r-1}l_1}}z_t\Ga^{l_1}_{lf_r}\psi_{t\iiv}.
\end{equation}

Finally,\\ $R^0 _{\iiv \i n+1;f_1;\cdots;f_r;n+1}=$
$$\begin{array}{rl} &\frac{\partial R^0 _{\iiv \i n+1;f_1;\cdots;f_r}}{\partial x^{n+1}}+
\Ga^0 _{n+1\,l_1}R^{l_1}_{\iiv \i n+1;f_1;\cdots;f_r} -\Ga^{l_1}_{n+1\,\iiv }R^{0 }_{l_1\i n+1;f_1;\cdots;f_r}-\Ga^{l_1}_{n+1\,\i }R^{0 }_{\iiv l_1n+1;f_1;\cdots;f_r}\\
&-\Ga^{l_1}_{n+1\,n+1}R^{0 }_{\iiv \i l_1;f_1;\cdots;f_r}
-\Ga^{l_1}_{n+1\,f_1}R^{0 }_{\iiv \i n+1;l_1;f_2:\cdots;f_r}-\cdots-\Ga^{l_1}_{n+1\,f_r}R^{0 }_{\iiv \i n+1;f_1:\cdots;f_{r-1};l_1}.\end{array}$$
Using this,  \eqref{LM2226},  claims 2) and 3), and the inductive hypothesis, we get
\begin{multline}\label{LM25} R^0_{\iiv \i n+1;f_1;\cdots;f_r;n+1}=\sum_{t\in T_{\i f_1...f_r}}\frac{\partial z_t}{\partial x^{n+1}}\psi_{t\iiv}
-\sum_{l_1=1}^{n}\sum_{t\in T_{ l_1f_1...f_r}}z_t\Ga^{l_1}_{n+1\i}\psi_{t\iiv}\\
-\sum_{l_1=1}^{n}\sum_{t\in T_{\i l_1f_2...f_r}}z_t\Ga^{l_1}_{n+1f_1}\psi_{t\iiv}-\cdots
-\sum_{l_1=1}^{n}\sum_{t\in T_{\i f_1...f_{r-1}l_1}}z_t\Ga^{l_1}_{n+1f_r}\psi_{t\iiv}.
\end{multline}

Combining \eqref{LM22} with \eqref{LM24} and \eqref{LM23} with \eqref{LM25} and using the fact that $\psi|_{\h'}=0$,
we see that claim 8) is true. The lemma is proved. $\Box$

From lemma \ref{LMlem2} it follows that $$\hol_0\subset\g^{4,\h,m,\psi}.$$ Thus, $$\hol_0=\g^{4,\h,m,\psi}.$$
The theorem is proved. $\Box$


\bibliographystyle{unsrt}

\end{document}